\documentclass[12pt, reqno]{amsart}

\usepackage[left=1.1in,top=1.2in,right=1.1in,bottom=1.2in]{geometry}
\usepackage{setspace,kantlipsum}
\usepackage{amsfonts}
\usepackage{amssymb}
\usepackage{amsmath}
\usepackage{amsthm}
\usepackage{latexsym}
\usepackage{graphicx}
\usepackage{xypic}
\usepackage{hyperref}
\usepackage{lineno}
\usepackage{tikz-cd} 
\usepackage{subcaption} 
\usepackage[bottom]{footmisc}
\renewcommand{\thefootnote}{\fnsymbol{footnote}}

\makeatletter
\@namedef{subjclassname@2020}{%
  \textup{2020} Mathematics Subject Classification}
\makeatother

\usepackage{doc}
\makeatletter
\def\maketitle{\par
      \begingroup
        \def\thefootnote{\fnsymbol{footnote}}%
      \setcounter{footnote}\z@
      \def\@makefnmark{\hbox to\z@{$\m@th^{\@thefnmark}$\hss}}%
      \long\def\@makefntext##1{\noindent
          \ifnum\c@footnote>\z@\relax
            \hbox to1.8em{\hss$\m@th^{\@thefnmark}$}##1%
          \else
          \hbox to1.8em{\hfill}%
            \parbox{\dimexpr\linewidth-1.8em}{\raggedright ##1}%
          \fi}
      \if@twocolumn\twocolumn[\@maketitle]%
      \else\newpage\global\@topnum\z@\@maketitle\fi
      \thispagestyle{titlepage}\@thanks\endgroup
      \setcounter{footnote}\z@
      \gdef\@date{\today}\gdef\@thanks{}%
      \gdef\@author{}\gdef\@title{}\gdef\@dedicatory{}}

\hypersetup{colorlinks,citecolor=red,pdfstartview=FitH, pdfpagemode=None}

\usepackage{xcolor, soul}

\numberwithin{equation}{section}

\theoremstyle{definition}
\newtheorem{theorem}{Theorem}[section]

\newtheorem{example}[theorem]{Example}
\newtheorem{corollary}[theorem]{Corollary}

\newtheorem{remark}[theorem]{Remark}

\renewcommand{\ge}{\geqslant}
\renewcommand{\le}{\leqslant}
\newcommand\norm[1]{\left\lVert#1\right\rVert}
\newcommand\ui{\| \!!!|}

\def\SO{{\rm SO}}
\def\SU{{\rm SU}}

\def\GL{{\rm GL}}
\def\SL{{\rm SL}}

\font\germ=eufm10

\def\a{\mbox{\germ a}}

\def\g{\mbox{\germ g}}

\def\k{\mbox{\germ k}}

\def\p{\mbox{\germ p}}

\def\ui{\|\hspace{-.25mm} |\,}

\def\R{\mathbb R}
\def\C{\mathbb C}
\def\Cnn{{\mathbb C}_{n\times n}}

\def\H{\mathbb H}
\def\N{\mathbb N}
\def\P{\mathbb P}
\def\U{\mathrm U}

\def\tr{{\rm tr\,}}

\def\ad{{\operatorname{ad}\,}}
\def\Ad{{\operatorname{Ad}\,}}

\def\conv{{\mbox{conv}\,}}

\begin{document}

\title[Inequalities and limits of weighted Spectral geometric mean]{Inequalities and limits of weighted spectral geometric mean}

\author{Luyining Gan}
\address{Department of Mathematics and Statistics\\ University of Nevada, Reno\\ Reno \\ NV 89557-0084\\ USA}
\email{lgan@unr.edu}

\author{Tin-Yau Tam}
\address{Department of Mathematics and Statistics\\ University of Nevada, Reno\\ Reno \\ NV 89557-0084\\ USA}
\email{ttam@unr.edu}

\dedicatory{In memory of Professor Boying Wang who passed away on August 20, 2021.
}
\thanks{Professor Boying Wang, a faculty member of Beijing Normal University during 1960--2002, was a pioneer of developing and promoting Multilinear Algebra in China.}



\begin{abstract}
We establish some new properties of spectral geometric mean.
In particular, we prove a log majorization relation between $\left(B^{ts/2}A^{(1-t)s}B^{ts/2} \right)^{1/s}$ and the $t$-spectral mean $A\natural_t B :=(A^{-1}\sharp B)^{t}A(A^{-1}\sharp B)^{t}$ of two positive semidefinite matrices $A$ and $B$, where $A\sharp B$ is the geometric mean, and the $t$-spectral mean is the dominant one. The limit involving $t$-spectral mean is also studied. 
We then extend all the results in the context of symmetric spaces of negative curvature.
\end{abstract}

\maketitle


\section{Introduction}
Let $\Cnn$ be the space of all $n \times n$ complex matrices, $\H_n$  the real space  of $n\times n$ Hermitian matrices,  $\P_n$  the set of $n\times n$ positive definite matrices in $\Cnn$
and $\U(n)$  the group of $n\times n$ unitary matrices. 
For any $X \in \Cnn$, both $e^X$ and $\exp X$  denote the exponential of $X$.
Given $A\in \Cnn$, we use $A\ge 0$ to denote that $A$ is positive semidefinite. Given $A, B\in \H_n$, denote by $A\le B$ the L\"owner order, that is, $B-A\ge 0$.  Given $X\in \Cnn$, denote by $\|X\|$ the spectral norm of $X$, that is, the largest singular value of $X$, and denote by $\sigma(X)$ the spectrum of $X$. If the eigenvalues of $X$ are all real, we write $\lambda(X) = (\lambda_1, \dots, \lambda_n)$, where $\lambda_1(X)\ge \cdots \ge \lambda_n(X)$ are the eigenvalues of $X$.

The metric geometric mean (geometric mean, for short) of $A, B\in \P_n$
\begin{equation}\label{eqn:gm}
A \sharp B := A^{1/2}(A^{-1/2}BA^{-1/2})^{1/2}A^{1/2},
\end{equation}
was first introduced by Pusz and Woronowicz \cite {PW75} in 1975 and further studied by Kubo and Ando~\cite{KA79} in the 1970s. Since then, it has been extensively studied. Though the definition \eqref{eqn:gm} looks awkward, it is indeed a natural generalization of the classical geometric mean $\sqrt {ab}$ of two positive numbers $a, b$ \cite{Bh07}.
Besides the algebraic formulation and properties, geometric mean has a rich geometric flavor which is due to the fact that $\P_n$ is a Riemannian manifold \cite {He78} and that the geometric mean $A \sharp B$  is the mid-point of  the unique geodesic joining $A$ and $B$ \cite{Bh07, LLT14}.

The spectral geometric mean (spectral mean, for short) of $A, B\in \P_n$ was introduced by Fiedler and Pt\'ak \cite {FP97} in {1997 and one of the formulations is}
\begin{equation}\label{eqn:sm}
A \natural B := (A^{-1}\sharp B)^{1/2}A(A^{-1}\sharp B)^{1/2}.
\end{equation}
They named it as spectral geometric mean because the square of $A \natural B$ is similar to $AB$, which means that the eigenvalues of their spectral mean are the positive square roots of the corresponding eigenvalues of $AB$ \cite [Theorem 3.2 and Remark 3.4] {FP97}.
As the spectral mean comes from the geometric mean, it possesses some important properties of the geometric mean and has been studied ~\cite{KL15, LL04, LMW21, Li12}. However, unlike the geometric mean, not many results have been obtained on the spectral mean. Thus, this paper aims to obtain new results
on the spectral mean and its extension, namely, the $t$-spectral mean. Some of the results are analogous to the  geometric mean.

For each $t \in [0, 1]$, the {\it $t$-metric geometric mean} ($t$-geometric mean, for short) and {\it $t$-spectral geometric mean} ($t$-spectral mean, for short) of $A$ and $B$ are naturally defined by
\begin{align}
A \sharp_t B & := A^{1/2}(A^{-1/2}BA^{-1/2})^{t}A^{1/2}, \quad {t \in [0, 1]},\label{t-sharp}\\
A \natural_t B & := (A^{-1}\sharp B)^{t}A(A^{-1}\sharp B)^{t}, \quad {t \in [0, 1]}. \label{t-spectral}
\end{align}
Lee and Lim~\cite{LL07} first introduced the $t$-spectral mean in 2007. In the same year, it was also studied by Ahn, Kim and Lim \cite [p.191]{AKL07} (also see \cite [p.446]{KL07}). Its further algebraic and geometric meaning has been recently studied by Kim~\cite{Ki21}.
When $t = 1/2$, they are abbreviated as $A \sharp_{1/2} B = A \sharp B$ and $A \natural_{1/2} B = A \natural B$.

Both~\eqref{t-sharp} and \eqref{t-spectral} are paths joining $A$ {(when $t=0$) and $B$ (when $t=1$)} in $\P_n$. Here is a good way to interpret \eqref{t-sharp}. Move the points $A$ and $B$ to 
$I$ and $A^{-1/2}BA^{-1/2}$ respectively, via the congruence action by $A^{-1/2}$:
\[
X\mapsto A^{-1/2}XA^{-1/2},\quad  X\in \P_n .
\]
The $t$-geometric mean of the commuting $I$ and $A^{-1/2}BA^{-1/2}$ is \[I\sharp_t (A^{-1/2}BA^{-1/2})= (A^{-1/2}BA^{-1/2})^{t}.\] Then apply the inverse action on $(A^{-1/2}BA^{-1/2})^{t}$
\[
Y\mapsto A^{1/2}YA^{1/2},\quad  Y\in \P_n,
\]
to have $A \sharp_t B = A^{1/2}(A^{-1/2}BA^{-1/2})^{t}A^{1/2}$. This nice feature follows from  the fact that $S\P_n:=\{A\in \P_n: \det A=1\}$ is a symmetric space \cite [p.208-209]{He78}, when it is identified with $\SL_n(\C)/\SU(n)$ via the polar decomposition, where $\SL_n(\C)$ is the special linear group over $\C$ and $\SU(n)$ is the special unitary group.

In this paper, we focus our study on $t$-spectral mean and organize the sections as follows. In Section 2 we provide a similarity property. 
In Section 3, we review a log majorization result of $t$-geometric mean and show that a similar result for the $t$-spectral mean. In particular, we prove that $\left(B^{ts/2}A^{(1-t)s}B^{ts/2} \right)^{1/s}$ is log majorized by $(A^s\natural_tB^s)^{1/s}$ for $s>0$, and  also log majorized by the $t$-spectral mean $A\natural_tB$ for a range of $s$ with respect to $t$.
In Section 4, we provide the limit of $t$-spectral mean when $p$ tends to $0$.
The results in Sections 2, 3 and 4 are then extended in the context of symmetric spaces associated with a noncompact semisimple Lie group in Section 5.
In Section 6, some remarks are given.

\section{Properties}

In this section, we establish a similarity property of $t$-spectral mean. 
Let $A, B\in \P_n$. We say that $A$ is {\it positively similar} to $B$ if there exists $C\in \P_n$ such that $A = CBC^{-1}$.
Let us first recall some basic properties of the $t$-spectral mean~\cite{Bh07, GLT21,LL07} in the following {theorem}.
\begin{theorem} \label{prop:t-spectral}
Let $A, B \in \P_n$ and $t \in [0, 1]$. Then
\begin{enumerate}
\item $(A \natural_t B)^{-1} = A^{-1} \natural_t B^{-1}$ and $A\natural_t B = B \natural_{1-t} A$.
\item $A^{-1} \sharp (A \natural_t B) = (B \natural_t A)^{-1} \sharp B = (A^{-1} \sharp B)^t$. \label{t-natural:2}
\item If $G_t = A^{-1} \sharp (A \natural_t B)$, then $A \natural_t B = G_t A G_t$ and $B \natural_t A ={ G_t^{-1}B G_t^{-1}}$.
\item $(A \natural_r B) \natural_t (A \natural_s B) = A \natural_{(1-t)r+ts} B$ for all ${t, r, s \in [0,1]}$.
\end{enumerate}
\end{theorem}

Fiedler and Pt\'ak \cite [Theorem 5.5(5)]{FP97} obtained a positively similarity relation between the geometric mean and the spectral mean.
\begin{theorem}{(Fiedler and Pt\'ak 1997)}
Given $A, B \in \P_n$, the geometric mean $A\sharp B$ is positively similar to $(A\natural B)^{1/2} U (A\natural B)^{1/2}$ for some $U\in \U(n)$.
\end{theorem}

Motivated by the above theorem, we prove that the geometric mean and $t$-spectral mean are positively similar, which is reduced to the result of Fiedler and Pt\'ak when $t=1/2$. 
\begin{theorem}\label{thm:similar}
Let $A, B \in \P_n$ and $t \in [0, 1]$. Then $A\sharp B$ is positively similar to \[(A\natural_{1-t} B)^{1/2} U (A\natural_t B)^{1/2}\] for some {$U\in \U(n)$}.
\end{theorem}

\begin{proof}
By Theorem \ref{prop:t-spectral}(3), 
$$A\natural_t B = G_tAG_t, \quad B\natural_t A = G_t^{-1}BG_t^{-1},$$ where $G_t = A^{-1} \sharp (A\natural_t B)\in \P_n$, that is, $$A= G_t^{-1} (A\natural_t B )G_t^{-1},\quad B = G_t (B\natural_t A) G_t.$$
Set $$W := G_t (B\natural_t A)^{1/2},\quad V := (A\natural_t B)^{-1/2}G_t.$$
Then \begin{equation}\label{eqn:WV}
 	W^* V^{-1} = (B \natural_t A)^{1/2}(A\natural_t B)^{1/2}.
 \end{equation}
Now 
\[
G_t (B\natural_t A) G_t = B = (A\sharp B) A^{-1}(A\sharp B)  = (A\sharp B) G_t (A\natural_t B)^{-1} G_t(A\sharp B),
\]
where the second equality holds because $A^{-1}\natural_t B = (A\sharp B)^t A^{-1} (A\sharp B)^t$, $0\le t\le 1$, is a curve joining $A^{-1}$ (when $t=0$) and $B$ (when $t=1$) in view of \eqref{t-spectral}. Then
\[
\begin{split}
&(A\natural_t B)^{-1/2}G_t^2(B\natural_t A) G_t^2 (A\natural_t B)^{-1/2}
\\=& (A\natural_t B)^{-1/2}G_t(A\sharp B) G_t (A\natural_t B)^{-1} G_t(A\sharp B)G_t (A\natural_t B)^{-1/2}
\\
 = &\left[(A\natural_t B)^{-1/2}G_t(A\sharp B) G_t (A\natural_t B)^{-1/2}\right]^2.
\end{split}
\]
Thus we have 
\[
\begin{split}
A\sharp B &= \left[(A\natural_t B)^{-1/2}G_t\right]^{-1} \left[(A\natural_t B)^{-1/2}G_t^2(B\natural_t A) G_t^2 (A\natural_t B)^{-1/2}\right]^{1/2}\left[G_t (A\natural_t B)^{-1/2}\right]^{-1} \\
&= V^{-1}(V W W^* V^*)^{1/2} (V^*)^{-1}.
\end{split}
\]
Set $R := VW$ and from~\eqref{eqn:WV} $R = VV^*(A\natural_t B)^{1/2}(B\natural_t A)^{1/2}$. The matrix $U := R^{-1} (R R^*)^{1/2}$ is unitary since
\[
UU^* = R^{-1} (R R^*)^{1/2}  (R R^*)^{1/2} (R^*)^{-1} = R^{-1} (R R^*)(R^*)^{-1} = I.
\]
 Then we have
\[
\begin{split}
(A\natural_{1-t} B)^{1/2} U (A\natural_t B)^{1/2} &= (B\natural_t A)^{1/2} U (A\natural_t B)^{1/2} \\
&= (B\natural_t A)^{1/2}  R^{-1} (R R^*)^{1/2} (A\natural_t B)^{1/2}\\
&= (B\natural_t A)^{1/2}\left[VV^*(A\natural_t B)^{1/2}(B\natural_t A)^{1/2}\right]^{-1}(VWW^*V^*)^{1/2}(A\natural_t B)^{1/2}\\
&= (A\natural_t B)^{-1/2}(V^*)^{-1}V^{-1}(VWW^*V^*)^{1/2}(V^*)^{-1}V^*(A\natural_t B)^{1/2}\\
&= \left[V^*(A\natural_t B)^{1/2}\right]^{-1}V^{-1}(VWW^*V^*)^{1/2}(V^*)^{-1}\left[V^*(A\natural_t B)^{1/2}\right]\\
&=  \left[V^*(A\natural_t B)^{1/2}\right]^{-1}(A\sharp B)\left[V^*(A\natural_t B)^{1/2}\right],
\end{split}
\]
where $V^*(A\natural_t B)^{1/2} = G_t^*\in \P_n$ since $V = (A\natural_t B)^{-1/2}G_t$. Thus we complete the proof.
\end{proof}

\section{Log Majorization}

Let $x = (x_1, x_2, \dots, x_n)$ and $y = (y_1, y_2, \dots, y_n)$ be in $\R^n$. Let $x^{\downarrow} = (x_{[1]}, x_{[2]}, \dots, x_{[n]})$ denote the rearrangement of the components of $x$ such that $x_{[1]} \ge  x_{[2]} \ge \cdots \ge x_{[n]}$. 
We say that $x$ is {\it majorized} by $y$ \cite{MOA11} , denoted by $x \prec y$,  if 
\[
\sum_{i=1}^k x_{[i]} \le \sum_{i=1}^k y_{[i]}, \quad k = 1, 2, \dots, n-1, \quad \text{and} \quad \sum_{i=1}^n x_{[i]} = \sum_{i=1}^n y_{[i]}.
\]
Among many equivalent conditions for majorization, the following is geometric in nature, noted by Rado \cite{R52}
 and
{A. Horn} \cite{Ho54}:
\[
x \prec y \quad \Leftrightarrow \quad \conv S_n \cdot x \subset \conv S_n \cdot y,
\]
where $\conv S_n \cdot x$ denotes the convex hull of the orbit of $x$ under the action of the symmetric group $S_n$. {See a good summary in Marshall, Olkin, and Arnold \cite[p.10-14, p.34] {MOA11}.}
When $x$ and $y$ are nonnegative vectors, we say that $x$ is {\it log majorized} by $y$, denoted by $x \prec_{\log} y$ if 
\[
\prod_{i=1}^k x_{[i]} \le \prod_{i=1}^k y_{[i]}, \quad k = 1, 2, \dots, n-1 \quad \text{and} \quad  \prod_{i=1}^n x_{[i]} = \prod_{i=1}^n y_{[i]}.
\]
{When} $x$ and $y$ are positive  vectors, $x \prec_{\log} y$ if and only if $\log x \prec \log y$, where $\log x:=(\log x_1, \log x_2, \dots, \log x_n)$. 

A natural way to extend the notion of log majorization from nonnegative real vectors  to positive semidefinite matrices is via their eigenvalues, that is,
given $X, Y\ge 0$, we write $X\prec_{\log} Y$ when $\lambda(X) \prec_{\log} \lambda (Y)$.  As a relation, log majorization is transitive, reflexive but not anti-symmetric, so it is not a partial order. Needless to say, it is  different from the L\"owner order $\le$ which is a partial order. We would like to point out that neither one implies the other. 

\begin{remark}\label{log-hyperbolic}
We would like to point out that $X\prec_{\log} Y$ can be extended to $X, Y$ which are diagonalizable with nonnegative eigenvalues. For example if $A, B$ are positive semidefinite, then $AB$ is diagonalizable with nonnegative eigenvalues, though $AB$ is not Hermitian in general.
\end{remark}

The $t$-geometric mean has been studied extensively and a lot of nice properties have been discovered. For example, the following result of Ando and Hiai \cite[Theorem 2.1]{AH94} gives a log majorization relation between the $t$-geometric mean of $r$-powers of positive semidefinite $A$ and $B$ and the $r$-power of the $t$-geometric mean of $A$ and $B$.  
\begin{theorem}(Ando and Hiai 1994)\label{Thm:AH}
For every $A, B \ge 0$ and $0\le t \le 1$,
\begin{equation}\label{eqn:geometric1}
A^r \sharp_t B^r  \prec_{\log} (A\sharp_t B)^r,  \quad r\ge 1,
\end{equation}
\begin{equation}\label{eqn:geometric2}
(A\sharp_t B)^r\prec_{\log}  A^r \sharp_t B^r,\quad 0< r\le 1,
\end{equation}
\begin{equation}\label{eqn:geometric3}
(A^p \sharp_t B^p)^{1/p}	\prec_{\log}  (A^q \sharp_t B^q)^{1/q}, \quad 0<q\le p.
\end{equation}
\end{theorem}

Theorem  \ref{Thm:AH} was extended to symmetric spaces of negative curvature by Liao, Liu and Tam \cite [Theorem 3.7]{LLT14}. See Remark \ref{geometry} for the geometry associated with \eqref{eqn:geometric1}, \eqref{eqn:geometric2} and \eqref{eqn:geometric3}. Motivated by Theorem \ref{Thm:AH}, we would like to know if analogous  relation holds for {the} $t$-spectral mean. The following theorem shows that such relation does exist, but in reverse order.
\begin{theorem}\label{Thm:spectral-power}\rm 
For every $A, B \ge 0$ and $0\le t \le 1$,
\begin{equation}\label{eqn:spectral1}
(A\natural_t B)^r \prec_{\log} A^r \natural_t B^r, \quad r\ge 1,
\end{equation}
\begin{equation}\label{eqn:spectral2}
A^r \natural_t B^r\prec_{\log} (A\natural_t B)^r ,\quad 0< r\le 1,
\end{equation}
\begin{equation}\label{eqn:spectral3}
 (A^q \natural_t B^q)^{1/q}\prec_{\log}(A^p \natural_t B^p)^{1/p}, \quad 0<q\le p,
 \end{equation}
that is, $p\to  (A^p \natural_t B^p)^{1/p}$ is a log majorization increasing function on $(0,\infty)$.

\end{theorem}

\begin{proof}
We first prove~\eqref{eqn:spectral1}.
We may consider $A, B>0$ by continuity argument.
It is easy to see
\[
\det ((A\natural_t B)^r) = (\det A)^{(1-t)r} (\det B)^{tr} = \det (A^r \natural_t B^r).
\]
Recall \cite[p.776-777]{MOA11}  that 
\begin{equation}\label{E:compound}
\prod_{i=1}^{k}\lambda_i(A)=\lambda_1(C_k(A)),\quad  k=1, \dots, n,
\end{equation}
 where $C_k(A)$ denotes the $k$th compound of $A \ge 0$. Thus, we need to show 
\[
\lambda_1(C_k((A\natural_t B)^r)) \le \lambda_1(C_k(A^r\natural_t B^r)), \quad k=1, \dots, n-1.
\]
Note that \cite [Lemma 1.2]{AH94},  \cite[p.781]{DAT16} $C_k(A\sharp_t B) = C_k(A)\sharp_t C_k(B)$ for $k=1, \dots, n$, and $C_k: \GL_n(\C) \to \GL_{n\choose k}(\C)$ is a group representation of the general linear group $\GL_n(\C)$. So from \eqref{t-spectral}, we have for $r\in \R$ and $t\in[0,1]$,
\begin{equation}\label{E:compound-sharp-r}
C_k(A^r\natural_t B^r) = C_k(A)^r\natural_t C_k(B)^r 
\end{equation}
and 
\[
C_k((A\natural_t B)^r) = (C_k(A)\natural_t C_k(B))^r.
\]
Hence it suffices to show that
\begin{equation}\label{eqn:lambda2}
\lambda_1(A\natural_t B)^r \le \lambda_1(A^r\natural_t B^r).
\end{equation}
By joint homogeneity of $t$-spectral mean, we have for $\alpha, \beta >0$
\begin{equation}
((\alpha A)\natural_t (\beta B))^r = \alpha^{r(1-t)}\beta^{rt}(A\natural_t  B)^r
\end{equation}
and
\begin{equation}
 ((\alpha A)^r\natural_t (\beta B)^r) = (\alpha^r A^r)\natural_t (\beta^r B^r) = \alpha^{r(1-t)} \beta^{rt} (A^r\natural_t B^r).
\end{equation}
In other words, both sides of~\eqref{eqn:lambda2} have the same order of homogeneity for $A, B$.
Thus to prove \eqref{eqn:lambda2}, we may show that 
\begin{equation}\label{E:lambda1}
\lambda_1(A^r\natural_t B^r)\le 1 \Rightarrow \lambda_1(A\natural_t B)^r \le 1.
\end{equation}
Suppose that  $\lambda_1(A^r\natural_t B^r)\le 1$, that is, $\lambda_1((A^{-r}\sharp B^r)^t A^r (A^{-r}\sharp B^r)^t)\le 1$, we have 
\[
(A^{-r}\sharp B^r)^t A^r (A^{-r}\sharp B^r)^t\le I
\]
and thus
\[
A^r \le (A^{-r}\sharp B^r)^{-2t}.
\]
For $r\ge 1$, that is, $0< 1/r\le 1$, we have
\begin{equation}\label{eqn:Asm}
A \le (A^{-r}\sharp B^r)^{-2t/r}.
\end{equation}
When $r\ge 1$, we have from \eqref{eqn:geometric1} 
\[
A^r \sharp_t B^r \prec_{\log} (A\sharp_t B)^r
\]
so that
\[
A^{-r} \sharp B^r \prec_{\log} (A^{-1}\sharp B)^r.
\]
As $r\ge 1$, from \eqref{eqn:Asm} and \eqref{eqn:geometric1} we have
\[
\lambda_1(A)\le \lambda_1 ((A^{-r} \sharp B^r)^{-2t/r})\le \lambda_1((A^{-1} \sharp B)^{-2t}),\]
which is equivalent to
\[
\lambda_1((A^{-1} \sharp B)^{2t})\lambda_1(A)\le 1.
\]
Since $A, (A^{-1} \sharp B)^{2t} \ge 0$, we have
\[
\lambda_1((A^{-1} \sharp B)^{2t})\lambda_1(A)
\ge \lambda_1((A^{-1} \sharp B)^{2t}A) 
= \lambda_1((A^{-1} \sharp B)^{t}A(A^{-1} \sharp B)^{t})
= \lambda_1(A\natural_t B).
\]
Thus we get $\lambda_1(A\natural_t B) \le 1$, that is, \eqref{E:lambda1} is established. Thus we complete the proof of \eqref{eqn:spectral1}. 
We omit the proofs of~\eqref{eqn:spectral2} and~\eqref{eqn:spectral3} due to the similar idea as that of~\eqref{eqn:spectral1}.
\end{proof}

Let us recall some interesting results in the following theorem. 
\begin{theorem} \label{Thm:AHA} \rm (Ando and Hiai 1994, Araki 1990)
Let $A, B \in \P_n$.
For any  $t \in [0, 1]$ and $s > 0$,
\begin{eqnarray}
A\sharp_{t}B &\prec_{\log} &e^{(1-t) \log A + t \log B} \label{E:Ando-Hiai}   \\
 &\prec_{\log}&\left(B^{ts/2}A^{(1-t)s}B^{ts/2} \right)^{1/s}.  \label{E:Araki}
\end{eqnarray}
\end{theorem}


The first inequality \eqref{E:Ando-Hiai} is a result of Ando and Hiai \cite[Corollary 2.3]{AH94} as the complementary counterpart of the famous Golden-Thompson inequality for Hermitian matrices $A$ and $B$: 
\[
\tr e^{A+B}\le \tr (e^Ae^B).
\] 
We remark that the complementary Golden-Thompson inequality 
\[
\tr (e^{pA}\sharp_{t}e^{pB})^{1/p} \le  \tr e^{(1-t) A + t  B},\quad p>0, \quad 0\le  t\le 1,
\]
was first proved by Hiai and Petz \cite {HP93} and then  extended to log majorization by Ando and Hiai \cite {AH94}.
The second inequality {\eqref{E:Araki}} follows from a result of Araki \cite{A90}.

Very recently Gan, Liu and Tam~\cite{GLT21}  have proved the following result which asserts that the $t$-geometric mean of two positive definite matrices is  log majorized by their $t$-spectral mean.

\begin{theorem} \label{Thm:spectral_geometric} \rm (Gan, Liu, and Tam 2021)
For all $A, B \in \P_n$ and $t \in [0, 1]$, 
\begin{equation}\label{spectral_geometric}
A \sharp_t B\prec_{\log} A \natural_t B.
\end{equation}
\end{theorem}

Motivated by Theorem \ref{Thm:AHA} and  Theorem \ref{Thm:spectral_geometric},
it is natural to ask whether an analogous log majorization relation exists between $\left(B^{ts/2}A^{(1-t)s}B^{ts/2} \right)^{1/s}$ and $A\natural_t B$, or between $e^{(1-t) \log A + t \log B}$ and $A\natural_t B$. The former would be a stronger result than the latter. We state this stronger result in Theorem~\ref{thm:spectral} in which the range of $s$ is specified. Before stating Theorem~\ref{thm:spectral}, we prove another interesting inequality between $\left(B^{ts/2}A^{(1-t)s}B^{ts/2} \right)^{1/s}$ and $(A^s\natural_t B^s)^{1/s}$ for all positive $s$.

\begin{theorem}\label{thm:spectral2}
Let $A, B \in \P_n$. For $t\in [0,1]$ and $s > 0$. We have
\begin{equation}\label{thm:natlog2}
\left(B^{ts/2}A^{(1-t)s}B^{ts/2} \right)^{1/s}\prec_{\log} (A^{s} \natural_t B^{s})^{1/s}.
\end{equation}
In particular, setting $s=1$ yields
\begin{equation}\label{chain}
A \sharp_t B\prec_{\log} e^{(1-t) \log A + t \log B}\prec_{\log} B^{t/2}A^{1-t} B^{t/2} \prec_{\log} A\natural_t B,\quad 0\le t\le 1.
\end{equation}
\end{theorem}

\begin{proof}
It is easy to see that for $s>0$, we have
\[
\det (A^{s} \natural_t B^{s})^{1/s} = (\det A)^{1-t} (\det B)^{t} =\det \left(B^{ts/2}A^{(1-t)s}B^{ts/2} \right)^{1/s}.
\]
Indeed, it is true for $s\in \R$.
By
\eqref{E:compound}, we need to show that for $s>0$,
\[
\lambda_1({C_k}((B^{ts/2}A^{(1-t)s}B^{ts/2} )^{1/s})) \le \lambda_1(C_k((A^{s} \natural_t B^{s})^{1/s})), \quad k=1, \dots, n-1.
\]
As the compound $C_k: \GL_n(\C) \to \GL_{n\choose k}(\C)$ is a group representation of the general linear group $\GL_n(\C)$, we have
\[
C_k((B^{ts/2}A^{(1-t)s}B^{ts/2})^{1/s}) = (C_k(B)^{ts/2}C_k(A)^{(1-t)s}C_k(B)^{ts/2})^{1/s}, 
\]
for all $k=1, \dots, n$.  From \eqref{E:compound-sharp-r}, we have
\[
C_k((A^{s} \natural_t B^{s})^{1/s}) = (C_k(A)^s\natural_t C_k(B)^s)^{1/s}.
\]
Hence it suffices to show that for all $A, B\in \P_n$ and $s>0$,
\begin{equation}\label{eqn:lambda3}
 \lambda_1((B^{ts/2}A^{(1-t)s}B^{ts/2})^{1/s}) \le \lambda_1((A^s\natural_t B^s)^{1/s}).
\end{equation}
Note that for $\alpha , \beta >0$,
\[
((\beta B)^{ts/2}(\alpha A)^{(1-t)s}(\beta B)^{ts/2})^{1/s} = \alpha^{1-t}\beta^t  (B^{ts/2}A^{(1-t)s}B^{ts/2})^{1/s}
\]
and
\[
 ((\alpha A)^s\natural_t (\beta B)^s)^{1/s} = ((\alpha^s A^s)\natural_t (\beta^s B^s))^{1/s} = \alpha^{(1-t)} \beta^{t} (A^s\natural_t B^s)^{1/s}.
\]
Thus $(B^{ts/2}A^{(1-t)s}B^{ts/2})^{1/s}$ and $A\natural_t  B$ have the same order of homogeneity for $A, B$. Then we may show that 
\begin{equation}\label{E:lambda1-natural}
\lambda_1((A^s\natural_t B^s)^{1/s}) \le 1 \Rightarrow \lambda_1((B^{ts/2}A^{(1-t)s}B^{ts/2})^{1/s})\le 1.
\end{equation}
Let $C(s):= A^{-s}\sharp B^s$.
Because $\lambda_1((A^s\natural_t B^s)^{1/s})\le 1$, which means that \[A^s\natural_t B^s = C^t(s)A^sC^t(s) \le I,\] 
we know \cite[p.114]{Bh97}
\begin{equation}
A^s\le C^{-2t}(s). \label{AC-2t}
\end{equation}
Applying the Riccati equation to $C(s):= A^{-s}\sharp B^s$, we have from \cite[p.11]{Bh07} and \eqref{AC-2t}
\begin{equation}
B^s =  C(s)A^sC(s) \le C^{2(1-t)}(s). \label{BCAC}
\end{equation}
Since $t\in [0, 1]$ and $1-t\in [0, 1]$, it follows from \eqref{AC-2t}, \eqref {BCAC}, and \cite[p.115]{Bh97} that
\[
A^{(1-t)s} \le C^{-2t(1-t)}(s) \le B^{-ts}.
\]
It amounts to $B^{ts/2}A^{(1-t)s}B^{ts/2} \le I$, which means that 
$ \lambda_1(B^{ts/2}A^{(1-t)s}B^{ts/2}) 
 \le 1$.  Thus we have
 \[
 \lambda_1((B^{ts/2}A^{(1-t)s}B^{ts/2})^{1/s}) = \lambda_1^{1/s}(B^{ts/2}A^{(1-t)s}B^{ts/2}) \le 1,
\]
that is, \eqref{E:lambda1-natural} is proved and thus the proof of~\eqref{thm:natlog} is completed.

Note that when $s=1$, it is natural to have
\[
B^{t/2}A^{1-t} B^{t/2} \prec_{\log} A\natural_t B.
\]
By Theorem \ref{Thm:AHA}, we have \eqref{chain} since $\prec_{\log}$ is transitive.
\end{proof}

\begin{theorem}\label{thm:spectral}
	Let $A, B \in \P_n$. For each chosen $t\in [0,1]$, let $0 < s \le \min\{1/t, 1/(1-t)\}$. We have
\begin{equation}\label{thm:natlog}
\left(B^{ts/2}A^{(1-t)s}B^{ts/2} \right)^{1/s}\prec_{\log} A \natural_t B.
\end{equation}
\end{theorem}

We omit the proof because it is almost identical to the proof of Theorem~\ref{thm:spectral2}. By replacing $C(s):=A^{-s}\sharp B^s$ by $C :=A^{-1}\sharp B$ in the above proof, \eqref{AC-2t} should be $A\le C^{-2t}$ and \eqref{BCAC} should be $B=CAC\le C^{2(1-t)}$. As $0<s\le \min\{1/t, 1/(1-t)\}$ for each chosen $t\in [0, 1]$, we know $ts, (1-t)s \in [0,1]$ and thus
\[
A^{(1-t)s} \le C^{-2st(1-t)} \le B^{-ts}.
\]
So we can derive Theorem~\ref{thm:spectral} that can compare $\left(B^{ts/2}A^{(1-t)s}B^{ts/2} \right)^{1/s}$ and $A\natural_t B$ in which the range of $s$ is specified.

\begin{remark} 
We remark that the  condition  $0<s \le \min\{1/t, 1/(1-t)\}$  for each chosen $t$ in Theorem \ref{thm:spectral} is more restrictive than the condition $0<s$ in Theorem \ref{Thm:AHA}. It is easy to see $\min\{1/t, 1/(1-t)\} \le 2$ for $t\in[0, 1]$. The following example shows that the upper bound $\min\{1/t, 1/(1-t)\}$ for $s$ is needed. Suppose that $t=1/2$. Then $\min\{1/t, 1/(1-t)\} = 2$. Now choose $s=2.1$ and
\[
A=
\begin{bmatrix}
563.2198 &  77.6893\\
   77.6893 &  71.7683
\end{bmatrix}
 \quad \text{and} \quad
B = 
 \begin{bmatrix}
40.7285 & -25.1376\\
  -25.1376  & 44.0770
  \end{bmatrix}.
\]
By MATLAB computation, we have
\[
\left(B^{ts/2}A^{(1-t)s}B^{ts/2} \right)^{1/s} =
\begin{bmatrix}
 135.6328 & -25.3588\\
  -25.3588  & 51.3716 
    \end{bmatrix}
\quad \text{and} \quad
A\natural_t B =
\begin{bmatrix}
139.2433 & -16.7122\\
  -16.7122 &  47.4272
\end{bmatrix},
\]
where the spectrum of $\left(B^{ts/2}A^{(1-t)s}B^{ts/2} \right)^{1/s}$ is $\{142.6760, 44.3285\}$ and the spectrum of $A\natural_t B$  is $\{142.1906, 44.4798\}$. Thus $\left(B^{ts/2}A^{(1-t)s}B^{ts/2} \right)^{1/s}\not\prec_{\log} A\natural_t B$.
\end{remark}

\begin{remark}\label{R:Araki}
For fixed $t\in [0,1]$,  Araki's result \cite {A90} and \cite [Theorem A]{AH94} asserts that
\[
\left(B^{tp/2}A^{(1-t)p}B^{tp/2} \right)^{1/p} \prec_{\log} \left(B^{tq/2}A^{(1-t)q}B^{tq/2} \right)^{1/q},\quad 0<p\le q.
\]
According to  \eqref {thm:natlog}, $A \natural_t B$ is an upper bound for the set of positive definite matrices $\left(B^{tp/2}A^{(1-t)p}B^{tp/2} \right)^{1/p}$ with respect to $\prec_{\log}$ for $0<p\le \min\{1/t, 1/(1-t)\}$. 
\end{remark}

From \eqref{E:Araki} and \eqref{thm:natlog2}, we have 
\begin{equation}\label{eqn:chain2}
	e^{(1-t)A + tB} \prec_{\log} \left(e^{tpB/2}e^{(1-t)pA}e^{tpB/2} \right)^{1/p}\prec_{\log} (e^{pA} \natural_t e^{pB})^{1/p},
\end{equation}
for $A, B\in \H_n$, $t\in [0, 1]$ and $p> 0$. As log majorization implies weak majorization (see \cite [p.42]{A94} and \cite [p.168]{MOA11}), we have the following corollary.

\begin{corollary}
If $A,B\in \H_n$ and $t\in [0,1]$, then for every $p> 0$
\[
\tr (e^{(1-t)A + tB}) \le \tr (e^{pA}\natural_t e^{pB})^{1/p}.
\]
Moreover, $\tr{(e^{pA}\natural_t e^{pA})^{1/p}}$ decreases to $\tr{(e^{(1-t)A + tB})}$ as $p\searrow 0$.
\end{corollary}

In the above corollary, $p\searrow 0$ means that $p\in \R$, as a variable, decreases to $0$.


\section{Limits of $t$-spectral mean}

Hiai and Petz \cite [Lemma 3.3] {HP93} determined the limit of $t$-geometric mean when $p$ tends to $0$:
\begin{equation}\label{E:limit-HP}
 \lim_{p \to 0} \left(e^{pA} \sharp_t e^{pB}\right)^{1/p}= e^{(1-t)A+tB}.
\end{equation}
Furthermore, $\ui \left(e^{pA} \sharp_t e^{pB}\right)^{1/p} \ui$ is increasing for any unitarily invariant norm $\ui\cdot\ui$, and thus 
$\ui \left(e^{pA} \sharp_t e^{pB}\right)^{1/p}\ui $ increases to $\ui e^{(1-t)A+tB} \ui$ as $p\searrow 0$. Ando and Hiai~\cite{AH94} proved that
\[
(e^{pA} \sharp_t e^{pB})^{1/p}\prec_{\log} e^{(1-t)A+tB},\quad p> 0.
\]
It means that $e^{(1-t)A+tB}$ is an upper bound for $(e^{pA} \sharp_t e^{pB})^{1/p}$ for all $p\ge 0$, with respect to $\prec_{\log}$. From Theorem \ref{Thm:AH}, we have 
\[
(A^p \sharp_t B^p)^{1/p}  \prec_{\log} (A^q \sharp_t B^q)^{1/q},\quad  0<q\le p.
\]
So we may write  
\begin{equation}\label{E:mono-down}
\left(e^{pA} \sharp_t e^{pB}\right)^{1/p}\nearrow_{\prec_{\log}} e^{(1-t)A+tB}, \quad \text{as}\quad p\searrow 0.
\end{equation}
Here the notation $\nearrow_{\prec_{\log}}$ means increasing with respect to $\prec_{\log}$.
On the other hand, recall that from \eqref{eqn:chain2}
\[
e^{(1-t)A+tB} \prec_{\log} (e^{pA} \natural_t e^{pB})^{1/p},\quad p> 0,
\]
that is, $e^{(1-t)A+tB}$ is a lower bound of $(e^{pA} \natural_t e^{pB})^{1/p}$ for all $p> 0$.
So it would be natural to 
ask 
if formulas similar to \eqref{E:limit-HP} and \eqref{E:mono-down} hold for $t$-spectral mean.
The answer is affirmative and is given in the following theorem. The limit in Theorem~\ref{Lem:spectral-exp} was proved in~\cite{AKL07} by differentiation method. Here, we provide another proof. 
Given two functions $f, g: \N \to \R^+$ from the set $\N$ of natural numbers to the set $\R^+$ of positive real numbers, the \emph{little-o} notation $f(n) = o(g(n))$ means that $g(n)$ grows much faster than $f(n)$ intuitively. Rigorously, it means that for all $\varepsilon>0$, there exists some $k \in \N$ such that $0\le f(n)<\varepsilon g(n)$ for all $n\ge k$. Thus $\lim_{n\to \infty}f(n)/g(n) = 0$.

\begin{theorem}\label{Lem:spectral-exp}
	If $A,B \in \H_n$  and $t\in [0, 1]$, then
\begin{equation}\label{exp-sharp}
	 \lim_{p \to 0} \left(e^{pA} \natural_t e^{pB}\right)^{1/p}=e^{(1-t)A+tB}.
\end{equation}
Moreover,
\begin{equation}\label{E:mono-up}\left(e^{pA} \natural_t e^{pB}\right)^{1/p}\searrow_{\prec_{\log}} e^{(1-t)A+tB}, \quad \text{as}\quad p\searrow 0.
\end{equation}
\end{theorem}

\begin{proof}
Suppose $p>0$. Let $w := -p<0$. Note that
\[\lim_{w \to 0^-} \left(e^{wA} \natural_t e^{wB} \right)^{1/w}
=\lim_{p \to 0^+} \left((e^{pA})^{-1} \natural_t (e^{pB})^{-1} \right)^{-1/p}
=  \lim_{p \to 0^+} \left(e^{pA} \natural_t e^{pB} \right)^{1/p},
\]
where the last equality follows from $(A\natural_t B)^{-1} = A^{-1}\natural_t B^{-1}$.
So it suffices to prove 
\[
e^{(1-t)A+tB} = \lim_{p \to 0^+} \left(e^{pA} \natural_t e^{pB}\right)^{1/p}.
\]
Now we consider $p \to 0^+$, for $p\in (0,1)$, write $p = (m+s)^{-1}$ and
\[
X(p) = e^{pA} \natural_t e^{pB} \quad \text{and} \quad Y(p) = e^{p[(1-t)A+tB]},
\]
where $m=m(p)\in \mathbb{N}$ and $s=s(p) \in [0,1)$.
Let $\|\cdot\|$ be the spectral norm, that is, $\|X\|$ is the largest singular value of $X\in \Cnn$. Then
\begin{equation}\label{norm_Y(p)}
\norm{Y(p)} \le e^{p[(1-t)\norm A+t\norm B]}.
\end{equation}
Since $o(p)/p \to 0$ as $p\to 0^+$, we have
\[
\begin{split}
e^{pA}\sharp_t e^{pB} &=  e^{pA/2}\left\{ \left[ \sum_{k=0}^{\infty} \frac{1}{k!}\left(-\frac{pA}{2}\right)^k\right]
\left[\sum_{k=0}^{\infty}\frac{(pB)^k}{k!}\right]
\left[\sum_{k=0}^{\infty}\frac{1}{k!}\left(-\frac{pA}{2}\right)^k\right]
\right\}^t e^{pA/2}\\
& = e^{pA/2}[I+p(B-A)+o(p)]^t e^{pA/2}\\
& = \left[I+\frac{pA}{2}+o(p)\right] [I+pt(B-A)+o(p)] \left[I+\frac{pA}{2}+o(p)\right]\\
& = I+p[(1-t)A+tB]+o(p).
\end{split}
\]
Thus 
\[
e^{pA}\sharp e^{pB} = I+\frac p2(A+B)+o(p).
\]
Applying it on $X(p)$ yields
\[
\begin{split}
X(p) &=  (e^{-pA}\sharp e^{pB})^t e^{pA} (e^{-pA}\sharp e^{pB})^t\\
 &=  [I+\frac p2(-A+B)+o(p)]^t
 [I+pA+o(p)] [I+\frac p2(-A+B)+o(p)]^t \\
  &=  [I+\frac {tp}2(-A+B)+o(p)]
 [I+pA+o(p)] [I+\frac {tp}2(-A+B)+o(p)].
\end{split}
\]
Since
$Y(p) =  e^{p[(1-t)A+tB]}$, we have
\begin{equation}\label{Ypm}
\norm{Y(p)^{1/p}-Y(p)^m} \to 0\quad  \text{when}\ p\to 0^+.
\end{equation}
 Write $Y(p) =  I+p[(1-t)A+tB] +o(p)$,
so we have 
\[
X(p)-Y(p) = o(p).
\]
By Theorem \ref{Thm:AHA}  (also see Kubo and Ando \cite{KA79}), we have
\[
\norm{e^{pA}\sharp_t e^{pB}} \le \|Y(p)\| \le e^{p[(1-t)\norm A+t\norm B]}.
\]
Then we obtain
\[
\begin{split}
\norm{X(p)} &= \norm{(e^{-pA}\sharp e^{pB})^t e^{pA} (e^{-pA}\sharp e^{pB})^t}\\
& \le \norm{e^{-pA}\sharp e^{pB}}^t \norm{e^{pA}} \norm{e^{-pA}\sharp e^{pB}}^t\\
& \le e^{p[(1-t)\norm A+t\norm B]} e^{p\norm{A}} e^{p[(1-t)\norm A+t\norm B]}\\
&= e^{p[(3-2t)\norm{A}+2t\norm{B}]}.
\end{split}
\]
Hence
\begin{eqnarray}
\norm{X(p)^{1/p}-X(p)^m}  &\le &\norm{X(p)}^m\norm{X(p)^s -I} \nonumber\\
&\le &e^{(3-2t)\norm{A}+2t\norm{B}}\norm{X(p)^s -I} \to 0\quad  \text{when}\ p\to 0^+\label{Xpm}
\end{eqnarray}
because $e^{(3-2t)\norm{A}+2t\norm{B}}$ is bounded and $X(p)\to I$ as $ p\to 0^+$, while $s\in[0,1)$.

We also have
\begin{eqnarray}
&&\norm{X(p)^m - Y(p)^m} \nonumber\\
&\le& m\norm{X(p)-Y(p)} \left(\max\{\norm{X(p)}, \norm{Y(p)}\}\right)^{m-1} \nonumber\\
&=& m\norm{X(p)-Y(p)} \norm{X(p)}^{m-1} \nonumber\\
&\le &(m+s)\norm{X(p)-Y(p)} e^{p(m-1)\xi} \quad \text{where} \ \xi :=(3-2t)\norm A + 2t\norm B\ge 0\nonumber\\
&= &\frac{1}{p} \norm{X(p)-Y(p)} e^{\frac {m-1}{m+s}\xi} \quad \text{as} \ p = \frac 1{m+s}\nonumber\\
&\le &\frac{1}{p} \norm{X(p)-Y(p)} e^\xi \to 0\quad  \text{when}\ p\to 0^+\label{XmYm}
\end{eqnarray}
because $X(p)-Y(p) =o(p)$. Hence by  \eqref{Ypm}, \eqref{Xpm}, and \eqref{XmYm}, we have
\begin{eqnarray*}
&& \norm{e^{pA}\natural_t e^{pB}-e^{(1-t)A+tB}}\\
& = & \norm{X(p)^{1/p}-Y(p)^{1/p}}\\
&\le& \norm{X(p)^{1/p}-X(p)^m}   +\norm{X(p)^m - Y(p)^m}+ \norm{Y(p)^{1/p}-Y(p)^m} \to  0 \quad \text{as}\ p\to 0^+.
\end{eqnarray*}
The proof of \eqref{exp-sharp} is completed. 


From \eqref{eqn:spectral3}, we have  $$(e^{pA} \natural_t e^{pB})^{1/p}\prec_{\log} (e^{qA} \natural_t e^{qB})^{1/q},\quad 0<p\le q,
$$
that is, $p\to  (e^{pA} \natural_t e^{pB})^{1/p}$ is a log majorization increasing function on $(0,\infty)$.
Together with \eqref{eqn:chain2}
\[
e^{(1-t)A+tB} \prec_{\log} (e^{pA} \natural_t e^{pB})^{1/p},\quad p> 0,
\]
\eqref{E:mono-up} is proved.
\end{proof}

\begin{corollary} 
For $p> 0$, $\ui \left(e^{pA} \natural_t e^{pB}\right)^{1/p} \ui$ decreases to $ \ui e^{(1-t)A+t B}\ui$ as $p \searrow 0$ for any unitarily invariant norm $\ui \cdot \ui$.
\end{corollary}

In view of Remark \ref{R:Araki} one may ask whether 
\[
\lim_{p\to 0} \left(B^{pt/2}A^{p(1-t)}B^{pt/2} \right)^{1/p}
\] exists or not when $A, B>0$, or equivalently,
\[
\lim_{p\to 0} (e^{ptB/2}e^{p(1-t)A}e^{ptB/2})^{1/p}
\] exists or not when $A, B\in \H_n$. The next theorem tells us that the answer is affirmative.
Note that 
\[
e^{ptB/2}e^{p(1-t)A}e^{ptB/2} = I+p[(1-t)A+tB]+o(p),
\]
 and
 \[
 \norm{e^{ptB/2}e^{p(1-t)A}e^{ptB/2}}\le \norm{e^{pA}\natural_t e^{pB}}.
 \]
 We can derive the following theorem by an almost identical proof of Theorem~\ref{Lem:spectral-exp} by replacing $ 
e^{pA} \natural_t e^{pB}$ with $e^{ptB/2}e^{p(1-t)A}e^{ptB/2}$ and Remark \ref{R:Araki}.

\begin{theorem}\label{thm:exp}
If $A,B \in \H_n$  and $t\in [0, 1]$, then
 \[
 \lim_{p \to 0} (e^{ptB/2}e^{p(1-t)A}e^{ptB/2})^{1/p} = e^{(1-t)A+t B}.
 \]
Moreover, for $p> 0$
\[
(e^{ptB/2}e^{p(1-t)A}e^{ptB/2})^{1/p}\searrow_{\prec_{\log}} e^{(1-t)A+tB}, \quad \text{as}\quad p\searrow 0.
\]
\end{theorem}

\begin{corollary}
For $p>0$, $\ui (e^{ptB/2}e^{p(1-t)A}e^{ptB/2})^{1/p}\ui$ decreases to $\ui e^{(1-t)A+t B}\ui $ as $p\searrow 0$ for any unitarily invariant norm $\ui \cdot \ui$.
\end{corollary}	

\begin{remark}
Recently Audenaert and Hiai  \cite {AH16} considered the convergence of the sequences $\{ (A^{p/2} B^p A^{p/2})^{1/p} \}_{p\in \N}$ and $\{(A^p\sharp B^p)^{2/q}\}_{p\in \N}$, where $A, B\ge 0$. They
proved that 
\[
\lim_{p\to \infty} (A^{p/2} B^p A^{p/2})^{1/p}
\]
exists but its explicit form  is  not known. They also showed that $$ \lim_{p\to \infty} (A^p\sharp B^p)^{2/q}$$ exists when $A, B$ are $2\times 2$ but the general case is unsettled. We do not know whether the sequence $\{\left(e^{pA} \natural_t e^{pB}\right)^{1/p}\}_{p\in \N}$ converges or not and it would be interesting to know the answer.
\end{remark}


\section{Kostant's Pre-order and Symmetric Spaces}

We first refer to~\cite{He78, Kn02} for the notation on symmetric spaces here. Let $G$ be a noncompact connected semisimple Lie group with Lie algebra $\g$. 
Let $\Theta$: $G\to G$ be a Cartan involution of $G$, and let $K$ be the fixed point set of $\Theta$, which is an analytic subgroup of $G$. Let $\theta = d\Theta$ be the differential map of $\Theta$. Then $\theta: \g \to \g$ is a Cartan involution and $\g = \k \oplus \p$ is a \emph{Cartan decomposition}, where $\k$ is the eigenspace of $\theta$ corresponding to the eigenvalue $1$ (and also the Lie algebra of $K$) and $\p$ is the eigenspace of $\theta$ corresponding to the eigenvalue $-1$ (and also an $\Ad K$-invariant subspace of $\g$ complementary to $\k$). The  Killing form $B$ on $\g$ is negative definite on $\k$ and positive definite  on $\p$, and the bilinear form $B_\theta$ defined by 
\[
B_{\theta} (X, Y) = -B(X, \theta Y), \quad X, Y \in\g
\]
is an inner product on $\g$. 
For each $X \in \g$, let $e^X = \exp X$ be the exponential of $X$. Let $P = \{e^X: \, X\in\p\}$. The map $\p \times K \to G$, defined by $(X, k) \mapsto e^Xk$, is a diffeomorphism.
So each $g \in G$ can be uniquely written as 
\begin{equation} \label{E:Cartan}
g  = pk 
\end{equation}
with $p = p(g) \in P$ and $k = k(g) \in K$. The decomposition $G = PK$ is called a \emph{Cartan decomposition of $G$}.

Let $*: G \to G$ be the diffeomorphism defined by $*(g) = g^* = \Theta (g^{-1})$. 
Because $K$ is the fixed point set of $\Theta$ and $\exp_{\g}: \p \to P$ is bijective,  we see that $p^*= p$ for all $p\in P$ and $k^*=k^{-1}$ for all $k\in K$. 
By the Cartan decomposition \eqref{E:Cartan}, we have for all $g \in G$
\begin{equation} \label{E:p}
p(g) = (gg^*)^{1/2}.
\end{equation}

An element $X \in \g$ is called real semisimple (resp., nilpotent) if $\ad X$ is diagonalizable over $\R$ (resp., nilpotent). An element $g \in G$ is called {\it hyperbolic} (resp., {\it unipotent}) if $g = \exp X$ for some real semisimple (resp., nilpotent) $X \in \g$; in either case $X$ is unique and we write $X = \log g$. An element $g \in G$ is called {\it elliptic} if $\Ad g$ is diagonalizable over $\C$ with eigenvalues of modulus $1$. According to \cite[Proposition 2.1]{Ko73}, each $g \in G$ can be uniquely written as
\begin{equation} \label{E:CMJD}
g = ehu,
\end{equation}
where $e$ is elliptic, $h$ is hyperbolic,  $u$ is unipotent, and the three elements $e, h$ and $u$ commute.
The decomposition \eqref{E:CMJD} is called the {\it complete multiplicative Jordan decomposition}, abbreviated as CMJD.

The Weyl group $W$ of $(\g, \a)$ acts simply transitively on $\a$ (and also on $A$ through the exponential map $\exp: \a \to A$).
For any real semisimple $X \in \g$, let $W(X)$ denote the set of elements in $\a$ that are conjugate to $X$, that is,
\[
W(X) = \Ad G(X) \cap \a.
\]
It is known from \cite[Proposition 2.4]{Ko73} that  $W(X)$ is a single $W$-orbit in $\a$.
Let conv\,$W(X)$ be the convex hull in $\a$ generated by $W(X)$. For each $g \in G$, define
\[
A(g) =\exp\mbox{conv}\, W(\log h(g)),
\]
where $h(g)$ is the hyperbolic component of $g$ in its CMJD.

The Kostant's pre-order $\prec_G$ on $G$  is defined (see \cite[p.426]{Ko73}) by setting  $\prec_G$ if 
\[
A(f) \subset A(g).
\]
This pre-order induces a partial order on the conjugacy classes of $G$. It is known from \cite[Theorem 3.1] {Ko73} that this pre-order $\prec_G$ does not  depend on the choice of $\a$. 

\begin{example}
If $G = \SL_n(\C)$, then $K = \SU(n)$, the special unitary group and $P$ is the space of positive definite matrices of determinant $1$. See \cite [p.430-431]{He78} for the CMJD of  $\SL_n(\C)$, which comes from the additive Jordan decomposition. The Kostant's pre-order   $\prec_ {\SL_n(\C)}$  is $\prec_{\log}$.
\end{example}

Let
\[ 
\pi: G\to G/K, \quad g\mapsto gK
\]
be the natural projection. Then $\p$ may be identified with the tangent space $T_o(G/K)$ of $G/K$ at the origin $o=eK$ via $d\pi$. Thus any $\Ad K$-invariant inner product $\langle\cdot,\cdot\rangle$ on $\p$ induces a unique $G$-invariant Riemannian metric on $G/K$ \cite [p.208-209]{He78}, that is, a Riemannian metric invariant under the natural action of $G$ on $G/K$ given by 
\[
(g,xK)\mapsto gxK.
\] Since $G$ is semisimple, the Killing form $B$ on $\g$ is nondegenerate. If $B$ is negative definite on $\k$ and positive definite  on $\p$, then the symmetric space $G/K$ is said to be of noncompact type. 

The map $G\to P$, $g\mapsto gg^*$, is onto. Because for any $g\in G$, it maps $gK$ to a single point $gg^*$, it follows that the map
\begin{equation}
\psi:\ G/K\ \to P, \quad gK \mapsto gg^* \label{mappsi}
\end{equation}
is a bijection. It is in fact a diffeomorphism by the Cartan decomposition $G=PK$. Via $\psi$, $P$ may
be identified with $G/K$, and so may be regarded as a symmetric space of noncompact type. Note that
for $p\in P$, $\psi^{-1}(p)=p^{1/2}K$, and $G$ acts on $P$ by 
\begin{equation}\label{G-action}
(g,p)\mapsto gpg^*.
\end{equation}

Let $G = PK$ be the Cartan decomposition of $G$. The map $p \mapsto p^{1/2}K$ {identifies} $P$ with $G/K$ as a symmetric space of noncompact type.  See \cite [p.349-350]{T16} for the example $G = \SL_n(\R)$ (the special linear group over $\R$), $K=\SO(n)$ (the special orthogonal group), $P$ is the set of real positive definite matrices of determinant $1$. 

The $t$-geometric mean of $p, q \in P$ was defined by Liao, Liu and Tam~\cite{LLT14}:
\begin{equation}\label{P:t-geometric mean}
p \sharp_{t}q = p^{1/2}\left(p^{-1/2}qp^{-1/2}\right)^{t}p^{1/2}, \quad 0 \le t \le 1.
\end{equation}
It is the unique geodesic in $P$ from $p$ (at $t=0$) to $q$ (at $t=1$). 
When $t = 1/2$, we abbreviate $p \sharp_{1/2} q$ as $p \sharp q$, 
Similarly, the $t$-spectral mean of $p, q \in P$ was defined by Gan, Liu, and Tam~\cite {GLT21}: 
\[
p \natural_t q = (p^{-1} \sharp q)^{t} p (p^{-1} \sharp q)^{t}, \quad 0 \le t \le 1.
\]
One may interpret  $p \natural_t q$ as the outcome of the action of $(p^{-1} \sharp q)^{t}$ on $p$ in view of the $G$-action \eqref{G-action} on $P$. 
Theorem \ref{Thm:AH} and Theorem \ref{Thm:AHA} were extended in  \cite [Theorem 3.6 and Theorem 3.5, respectively]{LLT14} to $P$:
\begin{theorem} (Liao, Liu, and Tam 2014)
Let $p,q \in P$ and $t \in [0, 1]$. Then 
\begin{eqnarray}
p^r\#_tq^r &\prec_G& (p\#_tq)^r, \quad r \ge 1,  \label{AH:P-r>1}\\
(p\#_tq)^r &\prec_G& p^r\#_tq^r, \quad 0 <  r \le 1,\label{AH:P-r<1} \\
(p^r\#_t q^r)^{1/r} &\prec_G & (p^s\#_t q^s)^{1/s},\quad 0<s\le r.\label{AH:P-mono}
\end{eqnarray}
\end{theorem}
The inequality \eqref{AH:P-mono} means that given $p,q\in P$, the function $s\mapsto (p^s\#_t q^s)^{1/s}$ is monotonic decreasing in the open interval $(0,\infty)$.

\begin{theorem} (Liao, Liu, and Tam 2014)
Let $p, q \in P$ and $t \in [0, 1]$. If $0 < r \le 1$, we have
\begin{eqnarray*}
p\#_{t}q &\prec_G  &e^{(1-t)  \log p + t \log q} \\
               &\prec_G  & (q^{tr/2}p^{(1-t)r}q^{tr/2})^{1/r}  
\end{eqnarray*}
\end{theorem}

Theorem \ref{Thm:spectral_geometric} was extended to $P$ \cite {GLT21}: 
\begin{theorem} \rm (Gan, Liu, and Tam 2021)
Let $p,q \in P$ and $t \in [0, 1]$. Then
\begin{equation}
p \sharp_t q\prec_G \ p \natural_t q.
\end{equation}
\end{theorem}

\begin{remark}\label{geometry} 
There is nice geometry hidden in \eqref{AH:P-r>1}, \eqref{AH:P-r<1}, and \eqref{AH:P-mono}. For instance, let us illustrate \eqref{AH:P-r>1}: $p^r\#_tq^r \prec_G (p\#_tq)^r$,  $r\ge 1$. From \eqref{P:t-geometric mean}, $p^r = e\sharp_r p$ and $q^r = e\sharp_r q$, where $e$ denotes the identity of $P$. So \eqref{AH:P-r>1} can be rewritten as
\begin{equation}\label{AH-hyperbolic}
(e\sharp_r p) \#_t (e\sharp_r q) \prec_G e \sharp_r (p\#_tq), \quad r\ge 1.
\end{equation}
Consider the geodesic triangle $\Delta(e, p, q)$ determined by the three points $e, p, q\in P$, abbreviated as $\Delta$.
Recall that $e\sharp_\mu p$, $\mu\in[0,1]$, is the geodesic emanating from $e$ with end point $p$ and clearly $e\sharp_r p$ is the point on the geodesic corresponding to time $\mu=r$. Similarly, $e\sharp_r q$ is the point on the geodesic $e\sharp_\mu q$, $\mu\in[0,1]$, emanating from $e$ with end point $q$, corresponding to time $\mu=r$. Now 
\begin{enumerate}
\item $p^r\#_tq^r=(e\sharp_r p) \#_t (e\sharp_r q)$ is the point on the geodesic $(e\sharp_r p) \#_\nu (e\sharp_r q)$, $\nu\in[0,1]$, emanating from $e\sharp_r p$  with end point $e\sharp_r q$, corresponding to time $\nu=t$. 
\item
$(p\#_tq)^r = e \sharp_r (p\#_tq)$ is the point on the geodesic $e \sharp_\xi (p\#_tq)$, $\xi\in[0,1]$, emanating from $e$  with end point $p\#_tq$, corresponding to time $\xi=r$.
\end{enumerate}
See the following figures. Figure \ref{fig:hyperbolic} illustrates that the Kostant's pre-order $\prec_G$ relating the different points $p^r\#_t q^r$ and  $(p\#_tq)^r$ in \eqref{AH:P-r>1}, or more precisely its equivalent form \eqref{AH-hyperbolic}, reflects the hyperbolic geometry of $P$. Figure~\ref{fig:euclidean} depicts a hypothetical Euclidean space in which
the two points $p^r\#_tq^r$ and $(p\#_tq)^r$ would be identical as the space is flat.


\begin{figure}[ht]
\includegraphics[width=10cm]{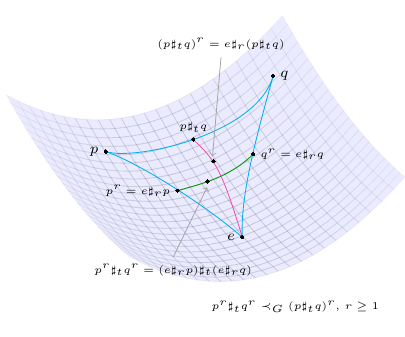}
\caption{\label{fig:hyperbolic} Hyperbolic space}
\end{figure}

\begin{figure}[ht]
\includegraphics[width=8cm]{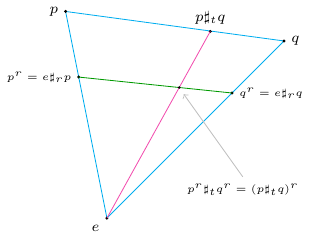}
\caption{\label{fig:euclidean} Euclidean space}
\end{figure}

\end{remark}

Using the technique in the proof of \cite[Theorem 3.5] {LLT14}, we can extend Theorems \ref{thm:similar}, \ref{Thm:spectral-power}, \ref{thm:spectral2}, \ref{thm:spectral}, \ref{Lem:spectral-exp}, and \ref{thm:exp} in the context of symmetric spaces $P$ of noncompact type as follows. 
\begin{theorem}
Let $p, q \in P$ and $t \in [0, 1]$. Then $p\sharp q$ is $P$-conjugate to $(q\natural_t p)^{1/2} {k_t} (p\natural_t q)^{1/2}$ for some $k_t\in K$. 
\end{theorem}

\begin{theorem}
For every $p, q \in P$ and $t \in [0, 1]$,
\[
(p\natural_t q)^r {\prec_G} \ p^r \natural_t q^r, \quad r\ge 1,
\]
or equivalently,
\[
p^r \natural_t q^r{\prec_G}\ (p\natural_t q)^r ,\quad 0< r\le 1,
\]
\[
(p^s \natural_t q^s)^{1/s}{\prec_G}\ (p^r \natural_t q^r)^{1/r}, \quad 0<s\le r.
\]
\end{theorem}

\begin{theorem}
Let $p, q \in P$. For $t \in [0, 1]$ and  $ s>0$,
\[
 \left(q^{ts/2}p^{(1-t)s}q^{ts/2} \right)^{1/s}{\prec_G} \ (p^s \natural_t q^s)^{1/s}.
\]
Moreover,
\[
p \sharp_t q\prec_G e^{(1-t) \log p + t \log q} \prec_G q^{t/2}p^{(1-t)}q^{t/2} \prec_G  p \natural_t q,\quad 0\le t\le 1.
\]
\end{theorem}

\begin{theorem}
Let $p, q \in P$. For each chosen $t\in [0,1]$, let $0 < s \le \min\{1/t, 1/(1-t)\}$. We have
\[
 \left(q^{ts/2}p^{(1-t)s}q^{ts/2} \right)^{1/s}{\prec_G} \ p \natural_t q.
\]
\end{theorem}

\begin{theorem}
For any $H,K \in \p $ and $t\in [0, 1]$,
\[	 
\lim_{p \to 0} (e^{pH} \natural_t e^{pK})^{1/p} = e^{(1-t)H+tK}.
\]
Moreover, 
\[(e^{pH} \natural_t e^{pK})^{1/p} \searrow_{\prec_G} e^{(1-t)H+tK} \quad \text{as} \quad p\searrow 0.\]
\end{theorem}

\begin{theorem}
For any $H,K \in \p $ and $t\in [0, 1]$,
 \[
 \lim_{p \to 0} (e^{ptK/2}e^{p(1-t)H}e^{ptK/2})^{1/p} = e^{(1-t)H+t K}.
 \]
Moreover,
\[(e^{ptK/2}e^{p(1-t)H}e^{ptK/2})^{1/p}\searrow_{\prec_G} e^{(1-t)H+tK}\quad \text{as} \quad p\searrow 0.\]
\end{theorem}


\section{Some remarks on L\"owner order and $t$-spectral mean}

The $t$-geometric mean $A \sharp_t B$ may be viewed as a function $\sharp_t: \P_n\times \P_n\to \P_n$ for any given $t\in[0,1]$:
$$(A,B) \mapsto A \sharp_t B=A^{1/2}(A^{-1/2}BA^{-1/2})^{t}A^{1/2}.$$ It is clearly a jointly continuous function.
Ando and Hiai~\cite [p.118]{AH94} mentioned that $t$-geometric mean is jointly (or coordinate-wisely) monotone with respect to the L\"owner order. 

\begin{theorem}[Ando and Hiai {1994}] 
\label{lem:geo_lowner} 
Let $A, B, C, D\in \P_n$. If $A\ge C, B\ge D$ and $t\in [0,1]$, then
\begin{equation}\label{Kim-LH}
A\sharp_t B \ge C\sharp_t D.
\end{equation}
\end{theorem}

An example was given by Dinh and Tam \cite [p.778]{DAT16} to show that the L{\"{o}}wner order in \eqref{lem:geo_lowner} cannot be replaced by $\prec_{\log}$.  The following example shows that the relation similar to Theorem \ref{lem:geo_lowner} is not true for the $t$-spectral mean.

\begin{remark}
Suppose $t=1/3$, $A=
\begin{bmatrix}
36.4987 & -34.0028\\
 -34.0028  & 39.8198
\end{bmatrix}$,
\[
B_1 = 
\begin{bmatrix}
6.8259  &-11.0027\\
-11.0027   &33.6773
\end{bmatrix} 
\quad \text{ and }
\quad
B_2 = 
\begin{bmatrix}
2.5166 &  -0.2222\\
-0.2222  &  3.4253	
\end{bmatrix}.
\]
It is easy to check $B_1\ge B_2$. Then
\[
A\natural_t B_1 = 
\begin{bmatrix}
21.5984 & -24.0515\\
  -24.0515  & 36.6270
\end{bmatrix} \quad \text{ and }
\quad
A\natural_t B_2 = 
\begin{bmatrix}
13.4040 & -10.9429\\
  -10.9429  & 15.7328	
\end{bmatrix},
\]
where the eigenvalues of $A\natural_t B_1-A\natural_t B_2$ are $-0.0213$ and $29.1098$.
Thus $A\natural_t B_1-A\natural_t B_2$ is not positive semidefinite.
\end{remark}

\begin{remark}
In 2012, Lim \cite{Li12} named the operator monotone property \eqref{Kim-LH} as L{\"{o}}wner-Heinz inequality. The readers should be alerted that in the literature  L{\"{o}}wner-Heinz inequality  is sometimes used (for example \cite [p.2-3]{Z02}) for the profound result of L{\"{o}}wner: 
\begin{equation}\label{LH-inequality}
A\ge B\ge 0 \Rightarrow A^r\ge B^r,\quad 0\le r\le 1.
\end{equation}  
 L{\"{o}}wner \cite{Lo34} first obtained the result in his 1934 seminal paper and Heinz \cite {H51} gave an alternative proof  in 1951. Since then, various  proofs of L{\"{o}}wner-Heinz inequality have been given by different authors. See the historical notes of Bhatia \cite[p.149-150] {Bh97}. That being said, it is true that \eqref{LH-inequality} follows from \eqref{Kim-LH}:
set $A=C=I$,   \eqref{Kim-LH} becomes $B^t = I\sharp_t B \ge  I\sharp_tD = D^t$ when $B\ge D$ and $0\le t\le 1$. In other words, L{\"{o}}wner-Heinz inequality can be written in the context of $r$-geometric mean:
\[
A\ge B\ge 0 \Rightarrow I\sharp_r A\ge I\sharp_r B,\quad 0\le r\le 1.
\]

\end{remark}

\begin{remark}

It is known that \cite[Theorem IX.2.6]{Bh97}
\begin{equation}\label{lambda1}
A, B\ge 0 \Rightarrow \lambda_1(A^sB^s) \le \lambda_1(AB)^s, \quad 0\le s\le 1.
\end{equation}
Both $AB$ and $A^sB^s$ are diagonalizable with  nonnegative eigenvalues so log majorization applies to them according to Remark \ref{log-hyperbolic}. Applying the compound matrix argument  on \eqref{lambda1}, we have
\begin{equation}\label{lambda1-2}
A, B\ge0 \Rightarrow A^sB^s \prec_{\log} (AB)^s, \quad 0\le s\le 1.
\end{equation}
Since $A^sB^s$ is  similar to $A^{s/2}B^sA^{s/2}$, \eqref{lambda1-2} can be rewritten as
\begin{equation}\label{lambda1-P}
A, B\ge 0 \Rightarrow A^{s/2}B^sA^{s/2} \prec_{\log} (A^{1/2}B^{1/2}A^{1/2})^s, \quad 0\le s\le 1.
\end{equation}
Let us confine ourselves in $\P_n$. The advantage of the form in \eqref{lambda1-P} is that all elements are now in $\P_n$. The group $\GL_n(\C)$ acts on $\P_n \cong \GL_n (\C)/\U(n)$ via \eqref{G-action}. So \eqref{lambda1-P} can be interpreted in the context of group action  on $\P_n$ and the $s$-geometric mean since  
\[
(A^{s/2},B^s)\mapsto  A^{s/2}B^sA^{s/2}\ \text{(group action)}
\]
and 
\[B^s = I\sharp_s B,\quad (A^{1/2}B^{1/2}A^{1/2})^s  = I\sharp_s A^{1/2}B^{1/2}A^{1/2}.
\]
Motivated by Figure~\ref{fig:hyperbolic}, the interested readers may draw a picture  to visualize the geometry associated with \eqref{lambda1-P}.

\end{remark}

\noindent{\bf Acknowledgement }
We are thankful to the anonymous referee for carefully reading of our paper and for giving constructive suggestions that helped us to greatly improve the paper.
\bibliographystyle{abbrv}
\bibliography{refs}

\end{document}